\documentclass[a4paper,12pt]{amsart}
\usepackage{amsfonts}
\usepackage{amssymb}
\usepackage{ifthen}
\usepackage{graphicx}
\nonstopmode \numberwithin{equation}{section}
\usepackage[usenames]{color}
\newtheorem{thm}{Theorem}[section]
\newtheorem{lem}{Lemma}[section]
\newtheorem{cor}{Corollary}[section]

\newtheorem{cl}{Claim}[section]
\newtheorem{ca}{Case}[section]
\newtheorem{sca}{Subcase}[section]
\newtheorem{scl}[section]{Subclaim}
\newtheorem{conj}[equation]{Conjecture}

\theoremstyle{definition}
\newtheorem{defn}{Definition}[section]
\newtheorem{op}[equation]{Open Problem}
\newtheorem{ques}[equation]{Question}
\newtheorem{rem}{Remark}[section]
\newtheorem{exam}[equation]{Example}

\newcounter {own}
\def\theown {\thesection       .\arabic{own}}

\newenvironment{pf}[1][]{%
 \vskip 3mm
 \noindent
 \ifthenelse{\equal{#1}{}}%
  {{\slshape Proof. }}%
  {{\slshape #1.} }%
 }%
{\qed\bigskip}

\newcounter{alphabet}
\newcounter{tmp}
\newenvironment{Thm}[1][]{\refstepcounter{alphabet}%
\bigskip%
\noindent%
{\bf Theorem \Alph{alphabet}}%
\ifthenelse{\equal{#1}{}}{}{ (#1)}%
{\bf .} \itshape}{\vskip 8pt}

\makeatletter
\newcommand{\Ref}[1]{\@ifundefined{r@#1}{}{\setcounter{tmp}{\ref{#1}}\Alph{tmp}}}
\makeatother

\newenvironment{Lem}[1][]{\refstepcounter{alphabet}%
\bigskip%
\noindent%
{\bf Lemma \Alph{alphabet}}%
{\bf .} \itshape}{\vskip 8pt}

\newcommand{\diam}{{\operatorname{diam}}}

\def\be{\begin{equation}}
\def\ee{\end{equation}}

\newcommand{\ben}{\begin{enumerate}}
\newcommand{\een}{\end{enumerate}}

\newcommand{\blem}{\begin{lem}}
\newcommand{\elem}{\end{lem}}
\newcommand{\bthm}{\begin{thm}}
\newcommand{\ethm}{\end{thm}}
\newcommand{\bcor}{\begin{cor}}
\newcommand{\ecor}{\end{cor}}
\newcommand{\beg}{\begin{exam}}
\newcommand{\eeg}{\end{exam}}
\newcommand{\begs}{\begin{examples}}
\newcommand{\eegs}{\end{examples}}
\newcommand{\bdefe}{\begin{defn}}
\newcommand{\edefe}{\end{defn}}
\newcommand{\bprob}{\begin{prob}}
\newcommand{\eprob}{\end{prob}}
\newcommand{\bques}{\begin{ques}}
\newcommand{\eques}{\end{ques}}
\newcommand{\bei}{\begin{itemize}}
\newcommand{\eei}{\end{itemize}}
\newcommand{\bcon}{\begin{conj}}
\newcommand{\econ}{\end{conj}}
\newcommand{\bop}{\begin{op}}
\newcommand{\eop}{\end{op}}

\newcommand{\bca}{\begin{ca}}
\newcommand{\eca}{\end{ca}}
\newcommand{\bsca}{\begin{sca}}
\newcommand{\esca}{\end{sca}}

\newcommand{\bcl}{\begin{cl}}
\newcommand{\ecl}{\end{cl}}

\newcommand{\bscl}{\begin{scl}}
\newcommand{\escl}{\end{scl}}

\newcommand{\bcons}{\begin{conjs}}
\newcommand{\econs}{\end{conjs}}
\newcommand{\bprop}{\begin{propo}}
\newcommand{\eprop}{\end{propo}}
\newcommand{\br}{\begin{rem}}
\newcommand{\er}{\end{rem}}
\newcommand{\brs}{\begin{rems}}
\newcommand{\ers}{\end{rems}}
\newcommand{\bo}{\begin{obser}}
\newcommand{\eo}{\end{obser}}
\newcommand{\bos}{\begin{obsers}}
\newcommand{\eos}{\end{obsers}}
\newcommand{\bpf}{\begin{pf}}
\newcommand{\epf}{\end{pf}}
\newcommand{\ba}{\begin{array}}
\newcommand{\ea}{\end{array}}
\newcommand{\beq}{\begin{eqnarray}}
\newcommand{\beqq}{\begin{eqnarray*}}
\newcommand{\eeq}{\end{eqnarray}}
\newcommand{\eeqq}{\end{eqnarray*}}

\newcounter{minutes}
\divide\time by 60
\newcounter{hours}
\multiply\time by 60 \addtocounter{minutes}{-\time}

\begin{document}

\bibliographystyle{amsplain}
\title{Deformations on symbolic Cantor sets and ultrametric spaces}

\author{Qingshan Zhou}
\address{Qingshan Zhou, School of Mathematics and Big Data, Foshan University,  Foshan, Guangdong 528000, People's Republic of China} \email{q476308142@qq.com}

\author{Xining Li}
\address{Xining Li,  Sun Yat-sen University, Department of Mathematics, Guangzhou 510275, People's Republic
of China} \email{lixining3@mail.sysu.edu.cn}

\author{Yaxiang Li${}^{~\mathbf{*}}$}
\address{Yaxiang Li,  Department of Mathematics, Hunan First Normal University, Changsha,
Hunan 410205, P.R.China}
\email{yaxiangli@163.com}

\def\thefootnote{}
\footnotetext{ \texttt{\tiny File:~\jobname .tex,
          printed: \number\year-\number\month-\number\day,
          \thehours.\ifnum\theminutes<10{0}\fi\theminutes}
} \makeatletter\def\thefootnote{\@arabic\c@footnote}\makeatother

\date{}
\subjclass[2010]{Primary: 30C65, 30F45; Secondary: 30C20} \keywords{Symbolic Cantor set, ultrametric space,  bilipschitz map, quasisymmetric map, quasim\"obius map.
\\
${}^{\mathbf{*}}$ Corresponding author}

\begin{abstract}  By introducing new deformations on symbolic Cantor sets and ultrametric spaces, we prove that doubling ultrametric spaces admit bilipschitz embedding into Cantor sets. If in addition the spaces are uniformly perfect, we show that they are quasisymmetrically equivalent to Cantor sets. Moreover, we provide a new proof for a recent work of Heer regarding quasim\"obius uniformization of Cantor set.
\end{abstract}

\thanks{The first author was supported by NNSF of
	China (No.11901090, 11571216), and by Department of Education of Guangdong Province, China (Grant nos.2018KQNCX285 and 2018KTSCX245). The second author   was   supported by NNSF of	China (No.2017TP1017). The third author was supported by  NNSF of
China (Nos.11601529,  11671127, 11971124)}

\maketitle{} \pagestyle{myheadings} \markboth{}{}

\section{Introduction and main results}\label{sec-1}
We start with the definition of {\it symbolic Cantor set} in \cite{DS} or \cite{Heer}. Let $F$ be a finite set with $k\geq 2$ elements and let $F^\infty$ denote the set of sequence $\{x_i\}_{i=1}^\infty$ with $x_i\in F$. Let $0<\lambda<1$. For two elements $x=\{x_i\}$, $y=\{y_i\}\in F^\infty$, we  define
$$L(x,y)=\sup\{I\in \mathbb{N}| \forall 1\leq i\leq I:x_i=y_i\}\;\;\;\;\mbox{and}\;\;\; \rho_\lambda(x,y)=\lambda^{L(x,y)}.$$
In particular, we have $L(x,x)=\infty$ and $L(x,y)=0$ if $x_1\neq y_1$. This defines an ultrametric on $F^\infty$. We call $(F^\infty,\rho_\lambda)$ the {\it symbolic $k$-Cantor set with parameter $\lambda$}. Here, a metric space $(X,d)$ is called {\it ultrametric} if for all $x,y,z\in X,$ we have the following {\it strong triangle inequality}:
\be\label{l-1} d(x,y)\leq \max\{d(x,z),d(z,y)\}.\ee

In an ultrametric space, any point of a ball is a center of the ball and all triangles are isosceles with at most one short side. It is well-known that the $p$-adic fields $\mathbf{Q_p}$ form a complete ultrametric space. There are many applications of ultrametric spaces in $p$-adic analysis, zeta function and fractal geometry; see \cite{DS,Ko,LF,LL,M,RK,S,WY}.

One easily observes that $F^\infty$ equipped with $\rho_\lambda$ is bounded and compact. It is natural to ask whether there is an unbounded ultrametric on $F^\infty$. We investigate this problem and get the following result.
\begin{thm}\label{thm1.1} Let $(F^\infty,\rho_\lambda)$ be the symbolic $k$-Cantor set with parameter $\lambda$. Then the function $\sigma_\lambda(x,y):=\lambda^{L(x,y)-L(x,o)-L(y,o)}$ on $F^\infty\setminus\{o\}$ is an unbounded ultrametric, where $o$ is a base point in $F^\infty$.
\end{thm}

\br  There are two deformations introduced in \cite{BHX} by Buckley, Herron and Xie.  The first class of flattening deformation on metric spaces  is a generalization of inversion on punctured $\mathbb{S}^n$. The second class of sphericalization deformation generalizes the conformal transformation from the Euclidean distance on $\mathbb{R}^n$ to the chordal distance on $\mathbb{S}^n$.  The original idea of these transformations follows from the work of Bonk and Kleiner \cite{BK02} in defining a metric on the one point compactification of an unbounded locally compact metric space.\er

In a recent work \cite{ZLL}, the authors proved that sphericalization and flattening sent quasi-metric spaces to quasi-metric spaces. It follows from Theorem \ref{thm1.1} that flattening also sent ultrametric spaces to ultrametric spaces. But the sphericalization deformation $s_p(x,y)=\frac{d(x,y)}{[1+d(x,p)][1+d(y,p)]}$  is not an ultrametric in general.
So it is natural to consider: for an unbounded ultrametric space, is there a bounded ultrametric on the space?

Motivated by this question, we introduce the following notation. Let $X$ be an ultrametric space with $a\in X$. We define the chordal metric $d_a$ on $\dot{X}=X\cup \{\infty\}$ by
\be\label{l-0}d_a(x,y)=\begin{cases}
\displaystyle\; \frac{d(x,y)}{\max\{1,d(x,a)\} \max\{1,d(y,a)\} },\;\;\;\;\;
\;\;\;\;\mbox{if}\;\;x,y\in X,\\
\displaystyle\;\;\;\;\;\frac{1}{\max\{1,d(x,a)\}},\;\;\;\;\;\;\;\;\;\;\;\;
\;\;\;\;\;\;\;\;\;\;\;\;\;\;\;\;\; \mbox{if}\;\;y=\infty\neq x,\\
\displaystyle\;\;\;\;\;\;\;\;\;\;0,\;\;\;\;\;\;\;\;\;
\;\;\;\;\;\;\;\;\;\;\;\;\;\;\;\;\;\;\;\;\;\;\;\;\;\;\;\;\;\;\;\;\;\;\;\;\;\mbox{if}\;\; x=\infty=y.
\end{cases}
\ee

We prove that the chordal metric $d_a$ is an ultrametric as follows.

\begin{thm}\label{thm1.2} Let $(X,d)$ be an ultrametric space with $a\in X$. Then the space $(\dot{X},d_a)$ is also an ultrametric space. \end{thm}

%

In \cite{BF}, Bonk and Foertsch  proved that a doubling compact ultrametric space admits a bilipschitz embedding into the symbolic $k$-Cantor set as follows.

\begin{Thm}\label{z-0}$($\cite[Proposition $6.3$]{BF}$)$
Suppose that $(X,d)$ is a doubling compact ultrametric space, then $(X,d)$ admits a bilipschitz embedding into the symbolic $k$-Cantor set $(F^\infty,\rho_\lambda)$ for sufficiently large $k$.
\end{Thm}

Note that a metric space has finite Assouad dimension (\cite{BF}) is equivalent to doubling. As an application of Theorems \ref{thm1.1} and \ref{thm1.2}, we obtain the following  unbounded analog for Theorem \Ref{z-0}.

\begin{thm}\label{thm1.3}
Suppose that $(X,d)$ is an unbounded, doubling and complete ultrametric space, then $(X,d)$ admits a bilipschitz embedding into the symbolic $k$-Cantor set $(F^\infty,\sigma_\lambda)$ for sufficiently large $k$, where $\sigma_\lambda(x,y):=\lambda^{L(x,y)-L(x,o)-L(y,o)}$ for some $o\in F^\infty$.
\end{thm}

\br We remark that  \cite[Corollary 2.8]{BDHM} is closely related to our Theorem \ref{thm1.3}. Without the doubling assumption for the ultrametric space, they proved that every separable ultrametric space is 3-bilipschitz embedding into a universal Cantor set.
\er

Moreover,  we consider the following quasisymmetric uniformization of symbolic Cantor set  which was proved by  David and Semmes in \cite{DS}.

\begin{Thm}\label{z-1}$($\cite[Proposition $15.11$]{DS}$)$
Every bounded, complete, doubling, uniformly perfect and uniformly disconnected metric space is quasi-symmetrically equivalent to the symbolic $2$-Cantor set $F^\infty$ equipped with the metric $\rho_\lambda$.
\end{Thm}

By Theorems \ref{thm1.1} and \ref{thm1.2}, we generalize Theorem \Ref{z-1} to the unbounded case.

\begin{thm}\label{thm1.4}
Every unbounded, complete, doubling, uniformly perfect and uniformly disconnected metric space is quasi-symmetrically equivalent to the symbolic $2$-Cantor set $(F^\infty\setminus\{o\},\sigma_\lambda)$, where $\sigma_\lambda(x,y):=\lambda^{L(x,y)-L(x,o)-L(y,o)}$ for some $o\in F^\infty$.
\end{thm}

Further, it should be mentioned that Heer recently generalized Theorem \Ref{z-1} and obtained the following quasim\"obius uniformization of Cantor set. As an application of Theorem \ref{thm1.4}, we provide a new proof for this result.

\begin{thm}\label{z-4}$($\cite[Theorem $5.2$]{Heer}$)$
Every complete, doubling, uniformly perfect and uniformly disconnected metric space is quasi-m\"{o}bius equivalent to the symbolic $2$-Cantor set $F^\infty$ equipped with the metric $\rho_\lambda$.
\end{thm}

The organization of this paper is as follows. In Section \ref{sec-2}, we recall some definitions and preliminary results. In Section \ref{sec-3}, we will prove our main results.

\section{Preliminaries}\label{sec-2}
\subsection{General Metric Space Information}Within the paper, we always  assume that $X$ is a metric space with a metric $d$. A metric space $X$ is called {\it doubling} if there exists an integer $n$ such that for all $r\leq \diam(X)$ and $x\in X$, there exist $n$ points $x_1,x_2,\dots,x_n$ with $B(x,r)\subset \cup_{i=1}^nB(x_i,r/2).$
We say that $X$ is {\it $C$-uniformly perfect}, if there exists a constant $C>1$ such that for each $x\in X$ and every $r>0$, $B(x,r)\setminus B(x, r/C)\not=\emptyset$ provided $X\setminus B(x,r)\not=\emptyset$.
$X$ is called {\it uniformly disconnected} if there exists a constant $\mu<1$ such that $X$ contains no $\mu$-chain, i.e. a sequence of (at least 3 distinct) points $(x_0,x_1,\ldots, x_n)$ such that $d(x_i,x_{i+1})\leq \mu d(x_0,x_n)$.

\subsection{Quasim\"obius, Quasisymmetric and Bilipschitz} Let $(X_1,d_1)$ and $(X_2,d_2)$ be two metric spaces. We say that $f$ is  {\it $L$-bilipschitz} if there exists $L\ge 1$ such that
$$d_1(x,y)/L\le d_2(f(x),f(y))\le L d_1(x,y).$$

Let $\eta:[0,\infty)\rightarrow [0,\infty)$ be a homeomorphism. We say that $f$ is  {\it $\eta$-quasisymmetric} if for $x,y,z\in X_0$,   we have
$$ \frac{d_2(f(x),f(z))}{d_2(f(x),f(y))}\leq \eta\left(\frac{d_1(x,z)}{d_1(x,y)}\right).$$

Given a metric space $(X,d),$ the {\it cross ratio} $r(x,y,z,w)$ of each four distinct points $x, y, z, w \in X$  is defined as
$$r(x,y,z,w)=\frac{d(x,z)d(y,w)}{d(x,y)d(z,w)}.$$

It is often convenient to consider cross ratios also in the extended space $\dot{X}$. If $x,y,z,w$ are points in $\dot{X}$ and if on the points $x,y,z,w$ is $\infty$, the cross ratio is defined by deleting the distances from $\infty$. For example $$r(x,y,z,\infty)=\frac{d(x,z)}{d(x,y)}.$$

Let $(X_1,d_1)$ and $(X_2,d_2)$ be two metric spaces, let $X_0\subset \dot{X}_1$, and let $f:(X_0,d_1)\rightarrow (\dot{X}_2,d_2)$ be a homeomorphism. Given a homeomorphism $\eta:[0,\infty)\rightarrow [0,\infty)$, we say that $f$ is {\it $\eta$-quasim\"obius} if for $x,y,z,w\in X_0,$ we have
$$r(f(x),f(y),f(z),f(w))\le \eta (r(x,y,z,w)).$$ If $f$ preserves all cross ratios, it is called a {\it M\"obius map}.  
For the properties of quasim\"obius  and quasisymmetric mappings see \cite{BK02,Vai,WZ17}. The following result concerning the relationship between quasim\"obius  and quasisymmetric mappings is very useful for our proofs. 

\begin{Thm}\label{Thm-l0} $($\cite[Theorem 3.10]{Vai}$)$
Suppose that $X$ is unbounded and that $f:X\to Y$ is $\theta$-quasim\"obius between two metric spaces. Then $f$ is $\theta$-quasisymmetric if and only if $f(x)\to\infty$ as $x\to \infty$.  If $X$ is any space
and if $f : \dot{X}\to \dot{Y}$ is $\theta$-quasim\"obius with $f(\infty)=\infty$, then $f|_X$ is $\theta$-quasisymmetric.
\end{Thm}

 Next, we also introduce some auxiliary results which will be used  later in our proofs.

\begin{lem}\label{newlem-1}
	The identity map $\psi: (\dot{X},d)\to (\dot{X},d_a)$ is M\"obius.
\end{lem}
\bpf
Given four points $x,y,z,w\in \dot{X}$, if one of them, let's say $x$,  is $\infty$, then we have $$\frac{d_a(\infty,z)d_a(y,w)}{d_a(\infty,y)d_a(z,w)}=\frac{d(y,w)}{d(z,w)}.$$ And if $x,y,z,w\in X$, then by the definition of $d_a$ we find that 
\begin{eqnarray*} & & \frac{d_a(x,z)d_a(y,w)}{d_a(x,y)d_a(z,w)}
\\ &=&\frac{d(x,z)}{\max\{1,d(x,a)\} \max\{1,d(z,a)\} }\frac{d(y,w)}{\max\{1,d(y,a)\} \max\{1,d(w,a)\} } 
\\ & & \cdot\frac{\max\{1,d(x,a)\} \max\{1,d(y,a)\} }{d(x,y)}\frac{\max\{1,d(z,a)\} \max\{1,d(w,a)\} }{d(z,w)}
\\ &=& \frac{d(x,z)d(y,w)}{d(x,y)d(z,w)}.
\end{eqnarray*} 
Hence the identity map $\psi$ is M\"obius.
\epf

\begin{lem}\label{newlem-2}
	The identity map $\varphi: (F^\infty,\rho_\lambda)\to (F^\infty\setminus \{o\} \cup\{\infty\},\sigma_\lambda)$
	is M\"obius with $\varphi(o)=\infty$, where $o$ is a base point in $(F^\infty,\rho_{\lambda})$.
\end{lem}
\bpf
Given four points $x,y,z,w\in F^\infty$, we consider two cases. If $x,y,z,w\in F^\infty\setminus\{o\}$, we note from the definition $$\sigma_\lambda(x,y)=\frac{\rho_{\lambda}(x,y)}{\rho_{\lambda}(x,o)\rho_{\lambda}(y,o)}$$ that
$$\frac{\sigma_\lambda(x,z)\sigma_\lambda(y,w)}{\sigma_\lambda(x,y)\sigma_\lambda(z,w)}
= \frac{\rho_{\lambda}(x,z)\rho_{\lambda}(y,w)}{\rho_{\lambda}(x,y)\rho_{\lambda}(z,w)}.
$$

If $x=o$ and $\varphi(o)=\infty$, then 
	$$\frac{\rho_\lambda(o,z)\rho_\lambda(y,w)}{\rho_\lambda(o,y)\rho_\lambda(z,w)}=\frac{\sigma_\lambda(y,w)}{\sigma_\lambda(z,w)}=\frac{\sigma_\lambda(\infty,z)\sigma_\lambda(y,w)}{\sigma_\lambda(\infty,y)\sigma_\lambda(z,w)}.$$
 
Hence the identity map $\varphi$ is M\"obius.
\epf

\begin{Lem}\label{Thm-z0} $($\cite[Lemma C]{WZ17}$)$
The image of a uniformly perfect metric space under a quasim\"obius map is also uniformly perfect.
\end{Lem}

\begin{Lem}\label{z-2}$($\cite[Proposition $15.7$]{DS}$)$
A metric space $(X,d)$ is uniformly disconnected if and only if there is an unltrametric $d'$ on $X$ that is bilipschitz equivalent to $d$.
\end{Lem}

\begin{Lem}\label{z-3}$($\cite[Theorem $1.1$]{Heer}$)$
Let $(X, d)$ be a doubling metric space and let $(Y,d')$ a metric space. Let $f:(X, d) \to (Y, d' )$ be a quasim\"obius homeomorphism. Then $(Y, d')$ is doubling.
\end{Lem}




\section{Proofs of  Theorem \ref{thm1.1} to Theorem \ref{z-4}}\label{sec-3}
\subsection{The proof of Theorem \ref{thm1.1}.} Since $(F^{\infty},\rho_{\lambda})$ is an ultrametric space, thus for any $x,y,z\in F^{\infty}\setminus\{o\}$ we have $$\rho_{\lambda}(x,y)\leq \max\{\rho_{\lambda}(x,z),\rho_{\lambda}(y,z)\},$$ which implies that
\be\label{l-1.1} \sigma_{\lambda}(x,y)=\frac{\rho_{\lambda}(x,y)}{\rho_{\lambda}(x,o)\rho_{\lambda}(o,y)}
\leq\max\{\frac{\rho_{\lambda}(x,z)}{\rho_{\lambda}(x,o)\rho_{\lambda}(o,y)},
\frac{\rho_{\lambda}(z,y)}{\rho_{\lambda}(x,o)\rho_{\lambda}(o,y)}\}.\ee

If $\rho_{\lambda}(z,o)\leq \min\{\rho_{\lambda}(x,o),\rho_{\lambda}(y,o)\}$, by (\ref{l-1.1}) we obtain the desired inequality
$$\sigma_{\lambda}(x,y)\leq \max\{\sigma_{\lambda}(x,z),\sigma_{\lambda}(z,y)\}.$$
It remains to consider the case $\rho_{\lambda}(z,o)> \min\{\rho_{\lambda}(x,o),\rho_{\lambda}(y,o)\}$. Without loss of generality, we may assume that $\min\{\rho_{\lambda}(x,o),\rho_{\lambda}(y,o)\}=\rho_{\lambda}(x,o)$. Then we have $$\rho_{\lambda}(z,o)>\rho_{\lambda}(x,o)$$ and so 
$$\rho_{\lambda}(z,o)=\rho_{\lambda}(x,z).$$ This deduces that $$\sigma_{\lambda}(x,z)=\frac{\rho_{\lambda}(x,z)}{\rho_{\lambda}(x,o)\rho_{\lambda}(z,o)}=\frac{1}{\rho_{\lambda}(x,o)}.$$
Moreover, since $\rho_\lambda(x,o)\leq \rho_\lambda(y,o)$, we get 
$$\rho_{\lambda}(x,y)\le \max\{\rho_{\lambda}(x,o),\rho_{\lambda}(y,o)\}=\rho_{\lambda}(y,o),$$
which shows   
$$\sigma_{\lambda}(x,y)=\frac{\rho_{\lambda}(x,y)}{\rho_{\lambda}(x,o)\rho_{\lambda}(y,o)}\le \frac{1}{\rho_{\lambda}(x,o)}=\sigma_{\lambda}(x,z).$$
Hence the proof of Theorem \ref{thm1.1} is complete.\qed


\subsection{The proof of Theorem \ref{thm1.2}.}
Fix $x,y,z\in X$. By symmetry, we only need to consider the following three cases.

\bca $\max\{d(x,a),d(y,a)\}\leq 1$.\eca
If $d(z,a)\leq 1$, then we have
$$d_a(x,y)=d(x,y)\leq \max\{d(x,z),d(y,z)\}=\max\{d_a(x,z),d_a(y,z)\}.$$

If $d(z,a)> 1$, then $$d_a(x,y)=d(x,y)\leq \max\{d(x,a),d(y,a)\}\leq 1=\max\{d_a(x,z),d_a(y,z)\}.$$

\bca  $d(x,a)\leq 1 < d(y,a)$.\eca
By the strong triangle inequality, we have $d(x,y)=d(y,a)$, which implies 
$$d_a(x,y)=\frac{d(x,y)}{\max\{1,d(x,a)\}\max\{1,d(y,a)\}}=\frac{d(x,y)}{d(y,a)}=1.$$

Next, we shall see that 
\be\label{l-2} \max\{d_a(x,z),d_a(y,z)\}=1.\ee
Indeed, if $d(z,a)\leq 1$, a similar argument as above gives that $d_a(y,z)=1$. Otherwise, if $d(z,a)>1$, then we also have $d_a(x,z)=1$. Therefore, in both cases we obtain the desired equality (\ref{l-2}).  

Hence we see that    
$$d_a(x,y)=1=\max\{d_a(x,z),d_a(y,z)\},$$
as required.

\bca  $\min\{d(x,a),d(y,a)\}>1$.\eca

If $d(z,a)\leq 1,$ thus we compute
\begin{eqnarray*}d_a(x,y)&=&\frac{d(x,y)}{d(x,a)d(y,a)}
\\&\leq& \max\{\frac{d(x,z)}{d(x,a)d(y,a)},\frac{d(z,y)}{d(x,a)d(y,a)}\}
\\&\leq& \max\{\frac{d(x,z)}{\max\{1,d(x,a)\}\max\{1,d(z,a)\}},\frac{d(z,y)}{\max\{1,d(y,a)\}\max\{1,d(z,a)\}}\}
\\&=& \max\{d_a(x,z),d_a(y,z)\}.\end{eqnarray*}

If $1<d(z,a)\leq \min\{d(x,a),d(y,a)\},$ then we find
\begin{eqnarray*}d_a(x,y)\leq \max\{\frac{d(x,z)}{d(x,a)d(y,a)},\frac{d(z,y)}{d(x,a)d(y,a)}\}
\leq\max\{d_a(x,z),d_a(y,z)\}.\end{eqnarray*}

If $d(z,a)> 1$ and $d(z,a)>\min\{d(x,a),d(y,a)\},$ by symmetry, we may assume that $\min\{d(x,a),d(y,a)\}=d(x,a)$. Then we have $d(z,a)\geq d(x,a)$ which implies $d(z,a)=d(x,z)$. Therefore, we obtain $$d_a(x,z)=\frac{d(x,z)}{d(x,a)d(z,a)}=\frac{1}{d(x,a)}.$$
On the other hand, since $d(x,a)\leq d(y,a)$, we get
$$d(x,y)\leq \max\{d(x,a),d(y,a)\}=d(y,a),$$
which yields \begin{eqnarray*}d_a(x,y)= \frac{d(x,y)}{d(x,a)d(y,a)}\leq \frac{1}{d(x,a)}=d_a(x,z).\end{eqnarray*}
Hence this proves Theorem \ref{thm1.2}.\qed

\subsection{The proof of Theorem \ref{thm1.3}.}
Consider the one-point extension space $\dot{X}=X\cup\{\infty\}$ equipped with the metric $d_a$ defined in (\ref{l-0}), where $a\in X$ is a base point. By Theorem \ref{thm1.2}, we see that $(\dot{X},d_a)$ is an ultrametric space. Then we know from Lemma \Ref{z-3} that $(\dot{X},d_a)$ is also doubling since from Lemma \ref{newlem-1} the identity map $(\dot{X},d)\to (\dot{X},d_a)$ is M\"obius. 

On the other hand, it is not difficult to see from \cite[Lemma 4.1.14]{HKST} that each complete doubling space is proper, i.e., every bounded closed ball is compact. Since $(\dot{X},d_a)$ is bounded, doubling and complete, thus we obtain that $(\dot{X},d_a)$ is compact.

Moreover, it follows from Theorem \Ref{z-0} that for some $k$ large enough, there exists an $L$-bilipschitz embedding $f$ from $(\dot{X},d_a)$ to $(F^\infty ,\rho_\lambda)$. Denote $f(\infty)=o$. We shall show that the induced map
$$f:(X,d)\to f(X)\subset (F^\infty ,\sigma_\lambda)$$
is $L^3$-bilipschitz. To this end, fix $x,y\in X$. By the definition of $d_a$ in (\ref{l-0}), we find that
$$d(x,y)=\frac{d_a(x,y)}{d_a(x,\infty)d_a(y,\infty)}$$
and so
\begin{eqnarray*} \sigma_\lambda(f(x),f(y))&=&\frac{\rho_\lambda(f(x),f(y))}{\rho_\lambda(f(x),f(\infty))\rho_\lambda(f(y),f(\infty))}
\\ &\leq& L^3 \frac{d_a(x,y)}{d_a(x,\infty)d_a(y,\infty)}
\\ &=& L^3 d(x,y).
\end{eqnarray*}
Similarly, we have 
$$\sigma_\lambda(f(x),f(y))\geq \frac{1}{L^3} d(x,y).$$
Hence this implies Theorem \ref{thm1.3}.\qed

\subsection{The proof of Theorem \ref{thm1.4}.}
Assume that $(X,\rho)$ is an unbounded, complete, doubling, uniformly perfect and uniformly disconnected metric space. Note first from Lemma \Ref{z-2} that there is an ultrametric metric $d$ on $X$ such that the identity map $\phi:(X,\rho)\to (X,d)$ is $L$-bilipschitz. Thus we see that $(X,d)$ is also unbounded, complete, doubling and uniformly perfect, because these properties are clearly bilipschitz invariant.

Next, we consider the one-point extension space $\dot{X}=X\cup\{\infty\}$ equipped with the metric $d_a$ defined in (\ref{l-0}), where $a\in X$ is a base point. A similar argument as the proof of Theorem \ref{thm1.3}, we see that $(\dot{X},d_a)$ is a compact ultrametric space. Because by Lemma \ref{newlem-1} we see that the identity map $\varphi:(X,d)\to (X,d_a)$ is M\"obius. 
Then it follows from Lemmas \Ref{Thm-z0} and \Ref{z-3} that $(X,d_a)$ is uniformly perfect and doubling as well.

Moreover, we claim that the metric completion $(\dot{X},d_a)$ of $(X,d_a)$ is also doubling and uniformly perfect. This can be seen as follows. 

On one hand, we observe from \cite[Lemma 4.1.14]{HKST} that $(\dot{X},d_a)$ is doubling with the same constant. On the other hand, assume that $(X,d_a)$ is $C$-uniformly perfect with constant $C>1$. It suffices to show that there is a constant $C_1>1$ depending only on $C$ such that, if $\dot{X}\setminus B_{d_a}(\infty,r)\neq \emptyset$ with $r>0$, then
$$B_{d_a}(\infty,r)\setminus B_{d_a}(\infty,r/C_1)\neq \emptyset,$$
where $B_{d_a}(\infty,r)\subset \dot{X}$ is the ball centered at $\infty$ with radius $r$ under the metric $d_a$. Since $(X,d)$ is unbounded, we choose a sequence of points $\{a_n\}\subset X$ with $1\leq d(a_n,a)\to\infty$ as $n\to\infty$. Thus we have 
$$d_a(a_n,\infty)=\frac{1}{\max\{1,d(a_n,a)\}}=\frac{1}{d(a_n,a)}\to 0,$$
as $n\to \infty$. This implies that the sequence $\{a_n\}$ is a $d_a$-Cauchy sequence and converges to $\infty$ in the space $(\dot{X},d_a)$. Since $\dot{X}\setminus B_{d_a}(\infty,r)\neq \emptyset$, there is some point $y\in X$ with 
$d_a(y,\infty)\geq r$. No loss of generality, we may assume that for all $n$
\[d_a(\infty,a_n)< \frac{r}{C}\] 
which implies that $d_a(a_n,y)=d_a(\infty,y)\geq r$.

Then, since $(X,d_a)$ is $C$-uniformly perfect, there is a point $z\in X$ such that
\[\frac{r}{C}\leq d_a(z,a_n)<r.\]
Therefore,  $$d_a(\infty,a_n)<\frac{r}{C}\leq d_a(\infty,z)=d_a(z,a_n)<r.$$ This yields that 
$$z\in B_{d_a}(\infty,r)\setminus B_{d_a}(z,\frac{r}{C})$$
and so $(\dot{X},d_a)$ is $C$-uniformly perfect.


Then it follows from Theorem \Ref{z-1} that there is an $\eta$-quasisymmetric homeomorphism $f$ from
$(\dot{X},d_a)$ onto the symbolic $2$-Cantor set $(F^\infty ,\rho_\lambda)$.
Denote $f(\infty)=o$ and
$$\sigma_\lambda(x,y)=\lambda^{L(x,y)-L(x,o)-L(y,o)}=\frac{\rho_{\lambda}(x,y)}{\rho_{\lambda}(x,o)\rho_{\lambda}(y,o)}$$
for all $x,y\in F^{\infty}\setminus \{o\}$. By Lemma \ref{newlem-2}, we observe  that the identity map
$$\psi:(F^\infty,\rho_\lambda)\to (F^\infty\setminus \{o\} \cup\{\infty\},\sigma_\lambda)$$
is M\"obius with $\psi(o)=\infty$. Consequently, we obtain an induced map
$$g:=\psi\circ f\circ\varphi:(\dot{X},d)\to (F^\infty\setminus \{o\} \cup\{\infty\},\sigma_\lambda),$$
which is quasim\"obius with $g(\infty)=\infty$.

Therefore, we observe from Theorem \Ref{Thm-l0} that $g_0=g|_X$ is actually quasisymmetric. Since the composition of bilipschitz and quasisymmetric maps is quasisymmetric, we find that the mapping 
$$g_0\circ \phi:(X,\rho)\to (F^\infty\setminus \{o\},\sigma_\lambda)$$
is quasisymmetric as well. 

Hence Theorem \ref{thm1.4} holds.\qed

\subsection{The proof of Theorem \ref{z-4}.} Assume that $(X,d)$ is a complete, doubling, uniformly perfect and uniformly disconnected metric space. If $X$ is bounded, then it follows immediately from Theorem \Ref{z-1} that $(X,d)$ is quasi-symmetrically equivalent to the symbolic $2$-Cantor set $F^\infty$ equipped with the metric $\rho_\lambda$.

It remains to assume that $X$ is unbounded, in this case we see from Theorem \ref{thm1.4} that $(X,d)$ is quasi-symmetrically equivalent to the symbolic $2$-Cantor set $(F^\infty\setminus\{o\},\sigma_\lambda)$, where $\sigma_\lambda(x,y):=\lambda^{L(x,y)-L(x,o)-L(y,o)}$ for some $o\in F^\infty$. Since from Lemma \ref{newlem-2} that the identity map
$(F^\infty,\sigma_\lambda)\to (F^\infty,\rho_\lambda)$ is M\"obius and since $o$ is the infinity point in $(F^\infty,\sigma_\lambda)$, we obtain that $(X,d)$ is quasim\"obius equivalent to $(F^\infty\setminus\{o\},\rho_\lambda)$ because the composition of quasisymmetric and M\"obius maps is quasim\"obius.

Hence the proof of Theorem \ref{z-4} is complete. \qed

\medskip

{\bf Acknowledgement.} The authors are indebted to the referee for the valuable suggestions.


\end{document}